\documentclass[12pt,leqno]{article}
\usepackage{amsmath,amssymb,amscd,latexsym,paralist}
\usepackage[all]{xy}

\usepackage[all]{xy}

\newcommand{\C}{{\mathbb{C}}}

\newcommand{\Q}{{\mathbb{Q}}}
\newcommand{\oQ}{\overline{\Q}}
\newcommand{\Z}{{\mathbb{Z}}}
\newcommand{\hZ}{\hat{\Z}}
\newcommand{\oZ}{\overline{\Z}}

\newcommand{\abb}{\mathrm{ab}}

\newcommand{\ann}{\mathrm{an}}

\newcommand{\id}{\mathrm{id}}

\newcommand{\Mor}{\mathrm{Mor}}
\newcommand{\rk}{\mathrm{rk}}
\newcommand{\Gal}{\mathrm{Gal}}
\newcommand{\GL}{\mathrm{GL}\,}
\newcommand{\Hom}{\mathrm{Hom}}
\newcommand{\uHom}{\underline{\Hom}}

\newcommand{\Ker}{\mathrm{Ker}\,}

\newcommand{\oN}{\overline{N}}
\newcommand{\tN}{\tilde{N}}
\newcommand{\Ob}{\mathrm{Ob}}

\newcommand{\uVec}{\mathbf{Vec}}

\newcommand{\Ch}{{\mathcal C}}
\newcommand{\Eh}{{\mathcal E}}
\newcommand{\Fh}{{\mathcal F}}

\newcommand{\Lh}{{\mathcal L}}
\newcommand{\Mh}{\mathcal{M}}
\newcommand{\Oh}{{\mathcal O}}

\newcommand{\eo}{\mathfrak{o}}
\newcommand{\eB}{\mathfrak{B}}

\newcommand{\eX}{{\mathfrak X}}
\newcommand{\eY}{{\mathcal Y}}

\newcommand{\eZ}{\mathcal{Z}}

\newcommand{\oF}{\overline{\mathbb{F}}}

\newcommand{\tC}{\tilde{C}}
\newcommand{\tD}{\tilde{D}}

\newcommand{\ohne}{\setminus}
\newcommand{\iso}{\stackrel{\sim}{\longrightarrow}}
\newcommand{\tei}{\, | \,}
\newcommand{\verk}{\mbox{\scriptsize $\,\circ\,$}}
\newcommand{\inv}{^{-1}}

\newtheorem{theorem}{Theorem}
\newtheorem{lemma}[theorem]{Lemma}
\newtheorem{prop}[theorem]{Proposition}
\newtheorem{defn}[theorem]{Definition}
\newtheorem{cor}[theorem]{Corollary}

\newtheorem{remark}[theorem]{Remark}
\newtheorem{conj}[theorem]{Conjecture}
\newtheorem{construction}[theorem]{Construction}

\newenvironment{rem}{\noindent {\bf Remark}}{}

\newenvironment{proof}{\noindent {\bf Proof}}{\mbox{}\hfill$\Box$}
\begin{document}
\title{Vector bundles on $p$-adic curves and parallel transport II}
\author{Christopher Deninger \and Annette Werner}
\date{\ }
\maketitle
\thispagestyle{empty}
\section{Introduction}
In \cite{dw1} we established a partial $p$-adic analogue of the Narasimhan--Seshadri correspondence between vector bundles and representations of the fundamental group. In \cite{fa} Faltings has even constructed a $p$-adic analogue of Simpson's theory of Higgs bundles. In his theory, $p$-adic representations of the fundamental group of a $p$-adic curve give rise to semistable $p$-adic Higgs bundles of slope zero. Faltings asks whether this condition is also sufficient for a $p$-adic Higgs bundle to come from a $p$-adic representation. 

Let $X$ be a smooth projective curve over $\oQ_p$. In \cite{dw1} we defined a category $\eB^{ps}_X$ of vector bundles on $X \otimes \C_p$ with ``potentially strongly semistable reduction of degree zero''. For these we constructed a theory of parallel transport along (homotopy classes of) \'etale paths and in particular representations of the fundamental group of $X$. Every vector bundle in $\eB^{ps}_X$ is semistable of slope zero on $X_{\C_p}$. Prompted by Faltings' question we wondered about the converse. In \cite{dw2} we showed that $\eB^{ps}_X$ equipped with the natural fibre functor in a point of $X_{\C_p}$ is a neutral Tannakian category just like the full category of semistable vector bundles of slope zero on $X_{\C_p}$.\\
Let $f : Y \to X$ be a possibly ramified finite covering of smooth projective curves over $\oQ_p$. It is well known that a vector bundle $E$ on $X_{\C_p}$ is semistable if and only if $f^* E$ is semistable on $Y_{\C_p}$. If $f$ is unramified then it is easy to see that $E$ is in $\eB^{ps}_X$ if and only if $f^* E$ is in $\eB^{ps}_Y$. However in the ramified case we do not know if this is true.

In the present paper we define for every $\mu \in \Q$ a category $\eB^{\mu}_X$ of vector bundles of slope $\mu$ on $X_{\C_p}$ with ``potentially strongly semistable reduction''. The bundles in $\eB^{\mu}_X$ are semistable of slope $\mu$ on $X_{\C_p}$ and $\eB^0_X$ contains $\eB^{ps}_X$. We show that a bundle $E$ is in $\eB^0_X$ if and only if $f^* E$ is in $\eB^0_Y$ even if $f$ is ramified, and similarly for non-zero slope. On the other hand we extend our theory of parallel transport from $\eB^{ps}_X$ to $\eB^0_X$. Since $\eB^0_X$ is again a neutral Tannakian category it may be reasonable to conjecture that every semistable vector bundle of slope zero on $X_{\C_p}$ lies in $\eB^0_X$ i.e. has potentially strongly semistable reduction.

The main difference in the definitions of $\eB^{ps}_X$ and $\eB^0_X$ lies in the meaning of the word ``potential''. For $E$ to be in $\eB^{ps}_X$ we allowed only pullbacks by finite \'etale coverings $\alpha$ of $X$ before $\alpha^* E$ was supposed to have ``strongly semistable reduction of degree zero''. For bundles $E$ in $\eB^0_X$ on the other hand we allow arbitrary ramified coverings $\alpha$ before $\alpha^* E$ is supposed to have ``strongly semistable reduction''. The problem in defining the parallel transport for $E$ in $\eB^0_X$ comes down to this: Using the theory of \cite{dw1} we can define a parallel transport between the fibres of $\alpha^* E$ for all \'etale paths on $Y_{\C_p}$. It is easy to see that this parallel transport descends to a parallel transport for $E$ along the \'etale paths in $U \subset X$, where $U$ is the complement of the ramification points of $\alpha$. In particular, for $x_0 \in U (\C_p)$ one obtains a representation $\rho : \pi_1 (U , x_0) \to \GL (E_{x_0})$. The main point is: This representation has no monodromy at the ramification points i.e. that it factors over $\pi_1 (X ,x_0)$. We were unable to prove this algebraically. Instead our proof uses Grothendieck's comparison theorem between algebraic and topological fundamental groups and some considerations on Riemann surfaces.

We were motivated by the following heuristic argument: A suitable small deformation of $\alpha$ and $Y$ would not affect the reduction of $\alpha^* E$ but replace the ramification locus of $\alpha$ on $X$ by a disjoint ramification locus. Thus we would have $U \cup U' = X$ for the new complement $U'$. By general arguments from \cite{dw1} and in particular propostion 34 the parallel transports on $U$ and $U'$ have to be compatible and glue to a parallel transport on $X$. In particular there would be no monodromy.

At the end of \cite{fa}, Faltings remarks that if the pullback with respect to a possibly ramified covering of one of his generalized representations $\rho$ of the fundamental group is a true representation, then $\rho$ is a true representation as well. This fact is compatible with a positive answer to his above question.

In the final section we extend another part of \cite{nase}, dealing with vector bundles of arbitrary slope, to the $p$-adic context. If $X$ has non-zero genus then for $r \ge 1$ there is a canonical central extension
\[
1 \longrightarrow \mu_r \longrightarrow \Gamma_r \longrightarrow \pi_1 (X , x_0) \longrightarrow 1
\]
where $\mu_r$ denotes the group of $r$-th roots of unity in $\C_p$. For any vector bundle $E$ of slope $\mu = d/r$ in $\eB^{\mu}_X$ we construct a representation
\[
\rho_{E, x_0} : \Gamma_r \longrightarrow \GL (E_{x_0}) \; .
\]
For this we use a similar topological argument as the one employed for proving the absence of monodromy above. It $E$ has degree zero then $\rho_{E,x_0}$ factors over $\pi_1 (X ,x_0)$ and is the representation defined previously for bundles in $\eB^0_X$.

The present paper originated from a question by Luis \'Alvarez-C\'onsul to the first author whether the Narasimhan--Seshadri correspondence for vector bundles of non-zero slopes had a $p$-adic analogue. We are very grateful to him for this insight. We would also like to thank Oscar Garc\'{\i}a-Prada very much for several explanations of the classical case.

\section{Vector bundles with strongly semistable reduction}
Recall that for a vector bundle $E$  on  a smooth, projective and connected curve $C$ over a field $k$ the slope is defined by  $\mu (E)
= {\deg (E)}/{\rk (E)}$. The bundle $E$ is called semistable
(respectively stable), if for all proper non-zero subbundles $F$ of $E$ the
inequality $\mu (F) \le \mu (E)$ (respectively $\mu (F) < \mu (E)$) holds.

If $\mbox{char} (k) = 0$, then pullback by any surjective morphism of
smooth connected projective curves preserves semistability of vector bundles. However in the case
$\mbox{char} (k) = p$, there exist vector bundles which are destabilized by the
Frobenius map. 

Assume that $\mbox{char} (k) = p$, and let
$F : C \to C$ be the absolute Frobenius morphism, defined by the $p$-power map
on the structure sheaf. 
In this case, a vector bundle $E$ on $C$ is called strongly semistable, 
if $F^{n*} E$ is semistable on $C$ for all $n
\ge 1$. 

Strong semistability  in characteristic $p$ is preserved by tensor products and by pullback with respect to $k$-morphisms between smooth connected projective curves, see \cite{mi}, \S\,5.

Now consider a one-dimensional proper scheme $Z$ over a field $k$ of characteristic $p$.
By $C_1 , \ldots , C_r$ we denote the irreducible components of $Z$
endowed with their reduced induced structures. Let $\tC_i$ be the normalization
of $C_i$, and write $\alpha_i : \tC_i \to C_i \to Z$ for the canonical morphism.
Note that the  $\tC_i$ are either points or smooth, projective and connected curves over $k$.

\begin{defn} \label{t1}
Let $E$ be a vector bundle on $Z$. Then $E$ is called strongly semistable, if one of the following equivalent conditions hold:
 
i) For all one-dimensional irreducible components $C_i$ of $Z$, the pullback $\alpha_i^* E$ is strongly semistable on $\tC_i$.

ii) For all  smooth, projective and connected curves $C$ over $k$ and all $k$-morphisms $\alpha: C \rightarrow Z$, the pullback 
$\alpha^* E$ is semistable on $C$.
\end{defn}

The equivalence of the two conditions can be shown as in \cite{dw2}, 12.2.4. Note however, that the slopes of the pullbacks
$\alpha_i^* E$ may depend on $i$. 

If $f: Z' \rightarrow Z$ is a $k$-morphism between one-dimensional proper $k$-schemes and $E$ is a strongly semistable
vector bundle on $Z$, then $f^* E$ is a strongly semistable vector bundle on $Z'$ (use ii) above). 

By $\oZ_p$ and $\eo$ we
denote the rings of integers in $\oQ_p$ and $\C_p$, and by $k = \oF_p$ the common residue field of
$\oZ_p$ and $\eo$. 

Let $X$ be a smooth, projective and connected curve over $\oQ_p$, and let
$X_{\C_p}$ be its base-change to $\C_p$. We call
any finitely presented, proper and flat scheme $\eX$ over $\oZ_p$ with generic fibre $X$ a model of $X$. Set $\eX_{\eo} = \eX \otimes_{\oZ_p} \eo$ and $\eX_k =
\eX \otimes_{\oZ_p} k$. 

\begin{defn}
  \label{t2} We say that a vector bundle $E$ of arbitrary slope on $X_{\C_p}$ has strongly
semistable reduction if $E$ is isomorphic to the generic fibre
of a vector bundle $\Eh$ on $\eX_{\eo}$ for some model $\eX$ of $X$, such that
the special fibre $\Eh_k$ is a strongly semistable vector bundle 
on the one-dimensional proper $k$-scheme $\eX_k$ in the sense of definition \ref{t1}.

We say that $E$ has potentially strongly semistable reduction, if there is a
finite  morphism $\alpha : Y \to X$ of connected smooth projective curves over
$\oQ_p$ such that $\alpha^*_{\C_p} E$ has strongly semistable reduction on $Y_{\C_p}$.

For all slopes $\mu \in \Q$ we denote by 
$\eB_X^\mu$ or simply $\eB^\mu$ the full subcategory of all vector bundles on $X_{\C_p}$ of slope $\mu$ which have potentially strongly semistable reduction.

\end{defn}

In \cite{dw1} and \cite{dw2} we considered the category $\eB_X^s$ of vector bundles on $X_{\C_p}$ with strongly semistable reduction of degree zero. Here $E$ has strongly semistable reduction of degree zero if there exists a model $\Eh$ on some $\eX_{\eo}$ as above, such that the pull-back of $\Eh_k$ to each normalized irreducible component of $\eX_k$ is strongly semistable of degree zero. In particular this implies that $E$ has degree zero. 
Note that a bundle of degree zero which has strongly semistable reduction in the more general sense of definition \ref{t2} may have different
degrees on the components of the special fibre. 

Besides we defined the category $\eB^{ps}_X$ as the category of vector bundles on $X_{\C_p}$ for which there exists a finite  \'etale covering $\beta: Y \rightarrow X$ by a connected curve $Y$ over $\oQ_p$ such that $\beta_{\C_p}^* E$ lies in $\eB_Y^s$. In the definition of $\eB^{\mu}$ on the other hand we allow ramified coverings. There are the following inclusions of categories
\[
\eB^s_X \subset \eB^{ps}_X \subset \eB^0_X \; .
\]

The following result provides a link between our old and our new categories. 

\begin{theorem}
\label{t3}

i) Let $E$ be a vector bundle of slope $\mu$ on $X_{\C_p}$. Then E has potentially strongly semistable reduction, i.e. $E$ lies in $\eB_X^\mu$, if and only if there exists a finite morphism $\alpha: Y \rightarrow X$ of smooth, projective, connected curves over $\oQ_p$ and a line bundle $L$ on $Y_{\C_p}$ such that the bundle $\alpha_{\C_p}^* E \otimes L$ lies in the category $\eB_Y^{s}$. 

ii) We can replace the category $\eB^s$ in i) by the category $\eB^{ps}$, i.e. $E$ lies in $\eB_X^\mu$ if and only if there exists a finite morphism $\alpha: Y \rightarrow X$ of smooth, projective, connected curves over $\oQ_p$ and a line bundle $L$ on $Y_{\C_p}$ such that $\alpha^*_{\C_p} E \otimes L$ lies in $\eB_Y^{ps}$. 
If $E$ satisfies this condition, then every line bundle $L$ on $Y_{\C_p}$ satisfying $\mbox{deg}\, L = - \mu \, \mbox{deg} \alpha$ has the property that $\alpha_{\C_p}^* E \otimes L$ lies in $\eB_Y^{ps}$.

iii) If $E_1, \ldots, E_r$ are finitely many vector bundles in $\eB_X^\mu$, then there is one covering $Y$ working for all of them. To be precise, there exists a finite dominant morphism $\alpha: Y \rightarrow X$ of smooth, projective, connected curves over $\oQ_p$ and line bundles $L_1,\ldots, L_r$ on $Y_{\C_p}$ such that the bundles $\alpha_{\C_p}^* E_i \otimes L_i$ lie in  $\eB_Y^{s}$. Moreover, for any line bundle $L$ on $Y_{\C_p}$ with $\deg L = -\mu \deg \alpha$ all bundles $\alpha^*_{\C_p} E_i \otimes L$ lie in $\eB^{ps}_Y$. 
\end{theorem}

\begin{proof} Let $E$ be a vector bundle of slope $\mu$ on $X_{\C_p}$. 

i) Let us first assume that there  is a finite covering $\alpha: Y \rightarrow X$ by a smooth, connected curve $Y$ and a line bundle $L$ on $Y_{\C_p}$ such the bundle $F = \alpha_{\C_p}^* E \otimes L$ lies in $\eB_Y^s$. This means that there exists  a model $\eY$ of $Y$ and a model $\Fh$ of $F$ on $\eY_{\eo} = \eY \otimes \eo$ such that the pullback of the special fibre
$\Fh_k$ to all normalized irreducible components $\tC_i$ of $\eY_k$ is strongly semistable of degree zero. By \cite{dw1}, theorem 5 there exists a model $\eY'$ of $Y$ dominating $\eY$ such that $L$ can be extended to a line bundle $\Lh$ on $\eY'_{\eo}$. If $\pi: \eY' \rightarrow \eY$ denotes the corresponding morphism, then the special fibre 
$(\pi^* \Fh)_k = \pi_k^* \Fh_k$ of $\pi^* \Fh$ is the pullback of a vector bundle which is strongly semistable of degree zero on all normalized irreducible components of $\eY$. Hence it is also semistable of degree zero on all normalized irreducible components of $\eY'$. 

Therefore we can substitute  $\eY$ by $\eY'$ and $\Fh$ by $\pi^* \Fh$ and assume that $L$ is the generic fibre of a line bundle $\Lh$ on $\eY_{\eo}$. 
The pullback of its special fibre $\Lh_k$ to any normalized irreducible component $\tC_i$ of $\eY_k$ is strongly semistable
since it is a line bundle. Hence the pullback of $(\Fh \otimes \Lh^{-1})_k = \Fh_k \otimes \Lh_k^{-1}$ to $\tC_i$ is strongly semistable. Since $\Fh \otimes \Lh^{-1}$ is a vector bundle on $\eY_{\eo}$ with generic fibre $\alpha_{\C_p}^* E$, we find that $\alpha_{\C_p}^* E$ has strongly semistable reduction in the sense of definition \ref{t2}. Therefore $E$ lies indeed in $\eB^\mu$. 

Now let us conversely assume that $E$ lies in $\eB_X^\mu$. We denote the degree of $E$ by $d$, and its rank by $r$. Hence 
$\mu = d/r$. By assumption, there exists a finite morphism $\alpha: Y \rightarrow X$,
a model $\eY$ of $Y$ and a vector bundle $\Fh$ on $\eY_{\eo}$ with generic fibre $\alpha_{\C_p}^*E$ such that its special fibre $\Fh_k$ is strongly semistable in the sense of definition \ref{t1}. Since $\eY$ is finitely presented over $\oZ_p$, it descends to a scheme over the ring of integers in a finite extension $K_0$ of $\Q_p$. Let $K$ be the maximal unramified extension of $K_0$, and let $\eo_K$ denote its ring of integers. It is a discrete valuation ring with residue field $k$. Then $\eY$ descends to a flat, proper and finitely presented $\eo_K$-scheme $\eY_{\eo_K}$. Its special fibre is $(\eY_{\eo_K})_k = \eY_k$. Note that $\eY_{\eo_K}$ is irreducible and reduced by \cite{liu}, Proposition 4.3.8. By lemma \ref{t4} proven below there exists an irreducible, semistable, proper $\eo_K$-scheme $\eZ$ with smooth generic fibre and a morphism of $\eo_K$-schemes
\[ \beta: \eZ \rightarrow \eY_{\eo_K} \]
with finite generic fibre $\beta_K$ such that for the irreducible components $D_1, \ldots, D_t$ of the special fibre $\eZ_k$ the following condition holds: If $\beta_k$ maps $D_i$ surjectively to an irreducible component $C_i$ of $\eY_k$ (i.e. not to a point), then the degree $k_i$ of the finite morphism 
\[{\beta_k}|_{D_i} : D_i \rightarrow  C_i,\]
where the components are endowed with their reduced structures, is a multiple of $r$. 
We denote by $\tD_i$ the normalization of $D_i$. If $\beta_k (D_i) = C_i$, then $\beta_k|_{D_i}: D_i \rightarrow C_i$ induces a morphism $\tD_i \rightarrow \tC_i$ of the same degree $k_i$.  

Now we look at the bundle $\beta_k^*  \Fh_k$. If $\beta_k$ maps the component $D_i$ of $\eZ_k$ to a point, then the pullback of $\beta_k^* \Fh_k$ to $D_i$ (and hence to $\tD_i$) is trivial, hence strongly semistable of degree zero.

On the other hand, if $\beta_k(D_i)$ is equal to the component $C_i$, then we denote by $(\beta_k^* \Fh_k)|_{\tD_i}$ the pullback of $\beta_k^* \Fh_k$ to $\tD_i$. The bundle $(\beta_k^* \Fh_k)|_{\tD_i}$ is equal to the pullback of $\Fh_k|_{\tC_i}$ via the map $\tD_i \rightarrow \tC_i$ induced by $\beta_k$. Since $\Fh_k|_{\tC_i}$ is strongly semistable, say of degree $e_i$, its pullback to $\tD_i$ is strongly semistable of degree $e_i k_i$. Hence the degree of $(\beta_k^* \Fh_k)|_{\tD_i}$ is a multiple of $r$. 

Now we want to construct for every irreducible component $D_i$ of $\eZ_k$ a line bundle $\Lh(i)$ on $\eZ$ such that the special fibre $\Lh(i)_k$ satisfies
\begin{eqnarray*}
\mbox{deg}(\Lh(i)_k|_{\tD_j}) = \delta_{ij},
\end{eqnarray*}
where $\delta_{ij}$ is the Kronecker delta. 

Since $\eZ$ is semistable, the singular locus $S$ in the special fibre $\eZ_k$ consists of finitely many closed points. 
Let $\eZ'= \eZ \backslash S$ be the complement. Then $\eZ'$ is a smooth $\eo_K$-scheme such that the special fibre $\eZ'_k$ is dense in $\eZ_k$. 
Hence every irreducible component $D_i$ of $\eZ_k$ contains a closed point $P_i \in \eZ'_k$. Since $k$ is algebraically closed, we have $P_i \in \eZ'_k(k)$. Note that the special fibre map
\[\eZ'(\eo_K) \rightarrow \eZ'_k(k)\]
is surjective, see e.g. \cite{blr}, 2.3, Proposition 5.
Hence we find some $\Delta_i \in \eZ'(\eo_K)$ with special fibre $P_i$.

 We identify $\Delta_i$ with its image in $\eZ'$,
which we regard as a Weil divisor. Since $\eZ'$ is regular, $\Delta_i$ gives rise to a line bundle
$\Oh(\Delta_i)$ on $\eZ'$. It  is trivial on the open subscheme $\eZ' \backslash \Delta_i$ of $\eZ$. Hence we can glue it with the trivial bundle on $\eZ \backslash \Delta_i$ and obtain a line bundle $\Lh(i)$ on $\eZ$. 
By construction we have
\[\mbox{deg}(\Lh(i)_k|_{D_j}) = 0 \mbox{, if } j \neq i \quad \mbox{and} \quad \mbox{deg}(\Lh(i)_k|{D_i}) = 1 \; .
\]
Therefore the line bundle $\Lh(i)$ on $\eZ$ satisfies indeed 
\begin{eqnarray*}
\mbox{deg}(\Lh(i)_k|_{\tD_j}) = \delta_{ij}.
\end{eqnarray*}

Now we consider again the vector bundle $\beta_k^* \Fh_k$ on the special fibre $\eZ_k$. We have shown that for all components $D_i$ the vector bundle $(\beta_k^* \Fh_k)|_{\tD_i}$ is  strongly semistable on $\tD_i$ such that $d_i = \mbox{deg} (\beta_k^* \Fh_k)|_{\tD_i}$ is a multiple of $r$. Let $\Lh^*$ denote the following line bundle on $\eZ$:
\[\Lh^* = \Lh(1)^{d_1 / r} \otimes \ldots \otimes \Lh(t)^{d_t/r}.\]
We consider the base changes $\eZ_{\oZ_p} = \eZ \otimes_{\eo_K} \oZ_p$ and $\eZ_{\eo} = \eZ \otimes_{\eo_K} \eo$, and denote by 
\[\beta_{\oZ_p} : \eZ_{\oZ_p} \rightarrow \eY_{\eo_K} \otimes_{\eo_K} \oZ_p = \eY\]
the base change map.  Then its generic fibre
\[\beta_{\oZ_p} \otimes_{\oZ_p} \oQ_p: Z = \eZ_{\oQ_p} \rightarrow \eY_{\oQ_p} = Y\]
is a finite morphism. Composing it with  $\alpha: Y \rightarrow X$, we get a finite morphism
\[ \gamma = \alpha \circ (\beta_{\oZ_p}\otimes_{\oZ_p}\oQ_p ) : Z \rightarrow X \]
 of smooth, projective, connected curves over $\oQ_p$. 
The vector bundle $\beta_{\eo}^* \Fh$ on $\eZ_{\eo}$ is a model of the vector bundle $\gamma_{\C_p}^* E$ on the generic fibre. Now let $\Lh = \Lh^* \otimes_{\eo_K} \eo$ be the base change of $\Lh^*$, and let $L = \Lh \otimes_{\eo} \C_p$ denote its generic fibre. Then the vector bundle $\beta_{\eo}^* \Fh \otimes \Lh^{-1}$ on $\eZ_{\eo}$ has generic fibre $\gamma_{\C_p}^* E \otimes L^{-1}$. Since the residue fields of $\eo_K$ and $\eo$ coincide, the special fibre of $\beta_{\eo}^* \Fh \otimes \Lh^{-1}$ is equal to $(\beta_k^* \Fh_k) \otimes \Lh^{* -1}_k$. Pulling back to the normalized irreducible component $\tD_i$ of $\eZ_k$ we get the vector bundle
\[(\beta_k^* \Fh_k)|_{\tD_i} \otimes \Lh^{* -1}_k|_{\tD_i}.\]
We know that $(\beta_k^* \Fh_k)|_{\tD_i}$ is strongly semistable of degree $d_i$ on $\tD_i$. As a line bundle, $\Lh^{* -1}_k|_{\tD_i}$ is also strongly semistable on $\tD_i$. By construction, its degree is equal to $-d_i / r$. Therefore
$(\beta_k^* \Fh_k)|_{\tD_i} \otimes \Lh^{* -1}_k|_{\tD_i}$ is strongly semistable of degree $0$. Hence
$\gamma_{\C_p}^* E$ lies indeed in $\eB_Z^s$. 

ii) Assume that there exists a finite covering $\alpha: Y \rightarrow X$ by a smooth, connected curve $Y$ and a line bundle $L$ on $Y_{\C_p}$ such the bundle $F = \alpha_{\C_p}^* E \otimes L$ lies in $\eB_Y^{ps}$. This means that there exists a finite \'etale covering $\beta: Y' \rightarrow Y$ by a projective, connected curve $Y'$ such that $\beta_{\C_p}^* (\alpha_{\C_p}^*E \otimes L)$ lies in $\eB_{Y'}^s$. Replacing $\alpha$ by $\alpha \circ \beta$ and $L$ by $\beta_{\C_p}^* L$, we are in the situation of i). Therefore $E$ lies in $\eB^\mu_X$. 

Now assume additionally that $M$ is any line bundle on $Y_{\C_p}$ satisfying $\deg M = -\mu \deg \alpha$. Then $L^{-1} \otimes M$ is a line bundle of degree zero on $Y_{\C_p}$. By \cite{dw1}, Theorem 12, $L^{-1} \otimes M$ lies in $\eB_Y^{ps}$. Since $\eB_Y^{ps}$ is stable under tensor product by \cite{dw1}, proposition 9, we find that $\alpha_{\C_p}^* E \otimes M$ lies indeed in $\eB_Y^{ps}$. 

iii) Let $E_1, \ldots, E_r$ be vector bundles in $\eB_X^\mu$. Then by i) there exist finite coverings $\alpha_i: Y_i \rightarrow X$ and line bundles $L_{i}$ on $Y_{i \C_p}$ such that $\alpha^*_{i \C_p} E \otimes L_i$ lies in $\eB_{Y_i}^s$ for all $i$. There is a finite covering $\alpha: Y \rightarrow X$ by a smooth, projective, connected curve $Y$ factoring over all $\alpha_i$. Since pullback via $Y \rightarrow Y_i$ maps $\eB_{Y_i}^s$ to $\eB_Y^s$ by \cite{dw1}, Proposition 9, the first claim follows. The second is a consequence of ii).
\end{proof} 

\begin{lemma}
\label{t4} Let $R$ be a discrete valuation ring with quotient field $K$ and residue field $k$, and
let $\eY$ be an  irreducible, reduced, flat and proper scheme  over $R$ with one-dimensional fibres, such that the generic fibre is smooth. Then for every integer $r \geq 1$  there exists an irreducible semistable proper  $R$-scheme $\eZ$  and an $R$-morphism $\beta: \eZ \rightarrow \eY$ such that the following conditions hold:

i) The generic fibre $\eZ_K$ is smooth, and $\beta_K: \eZ_K \rightarrow \eY_K$ is finite.

ii) For any irreducible component $D$ of $\eZ_k$ such that the special fibre $\beta_k$ maps $D$ surjectively to an irreducible component $C$ of $\eY_k$ (i.e. not to a point), the degree of the corresponding map of irreducible curves (endowed with their reduced structures) $\beta_k|_D: D \rightarrow C$ is a multiple of $r$. 

\end{lemma}

\begin{proof} Let $\eY'$ be the normalization of $\eY$  in its function field. Then the corresponding morphism $\pi: \eY' \rightarrow \eY$ is finite and an isomorphism on the generic fibre. Since it suffices to show our claim for the normalization $\eY'$, we can assume that $\eY$ is normal. Let $C_1, \ldots, C_t$ be the irreducible components of the special fibre $\eY_k$, and denote by $\eta_i$ their generic points. For any $i$ the local ring $\Oh_{\eY, \eta_i}$ is a discrete valuation ring in the function field $K(\eY)$. The residue field of $\Oh_{\eY, \eta_i}$ is also the residue field of the local ring $\Oh_{\eY_k, \eta_i}$ in the special fibre. Now $\Oh_{\eY_k, \eta_i}$ is a noetherian local ring of dimension $0$, hence it is Artinian and its residue field is simply the corresponding reduced ring $\Oh_{\eY_k, \eta_i}^{red}$. If we endow the component $C_i$ with its induced reduced structure, we have $\Oh_{C_i, \eta_i} = \Oh_{\eY_k, \eta_i}^{red}$. Hence the residue field of $\Oh_{\eY, \eta_i}$ is a function field of dimension one over $k$. Therefore it is the quotient field of a factorial ring (even of a discrete valuation ring). In particular, there exist irreducible polynomials, e.g. Eisenstein polynomials, of arbitrary degree over the residue field of $\Oh_{\eY, \eta_i}$. Hence there exists a monic polynomial $f(T) \in \Oh_{\eY, \eta_i}[T]$ of degree $r$ such that its reduction modulo the valuation ideal is irreducible over the residue field.

Then 
\[B_i = \Oh_{\eY, \eta_i} [T] / (f(T))\]
is a discrete valuation ring dominating $\Oh_{\eY, \eta_i}$ such that the corresponding extension of residue fields has degree $r$, see \cite{s} Ch. I, Proposition 15. Let $L_i$ be its quotient field. Then $L_i$ is an extension of $K(\eY)$ of degree $r$, such that $B_i$ is the only discrete valuation ring in $L_i$ dominating $\Oh_{\eY, \eta_i}$.

Let $L$ be the compositum of the extension fields $L_i$ in some algebraic closure of $K(\eY)$. Let $\eZ$ be the integral closure of the model $\eY$ in $L$. Then there is a finite morphism $\beta: \eZ \rightarrow \eY$.
If $D$ is an irreducible component of the special fibre $\eZ_k$ with generic point $\mu$,  then $\beta_k (D)$ is an irreducible component $C_i$ of $\eY_k$. Hence $\Oh_{\eZ, \mu}$ is a discrete valuation ring in $L$ dominating $\Oh_{\eY, \eta_i}$. Its intersection with $L_i$ is a discrete valuation ring in $L_i$ dominating $\Oh_{\eY, \eta_i}$, hence it is equal to $B_i$. Therefore the degree of the extension of residue fields corresponding to $\Oh_{\eY, \eta_i} \subset \Oh_{\eZ, \mu}$ is a multiple of $r$. 

If we endow the component $D$ of $\eZ_k$ with its induced reduced structure, the function field $K(D)$ of the curve $D$ is equal to the local ring $\Oh_{D, \mu}$. We have seen above that $\Oh_{D, \mu}$ coincides with the reduced ring associated to $\Oh_{\eZ_k, \mu} = \Oh_{\eZ,\mu} \otimes_R k$.  Hence the residue field of $\Oh_{\eZ, \mu}$ is equal to $\Oh_{D, \mu} = K(D)$.

Therefore the degree of the map $\beta_k|_{D}: D \rightarrow C_i$ is just the residue degree of $\Oh_{\eY, \eta_i} \subset \Oh_{\eZ, \mu}$ and hence a multiple of $r$. 

To conclude the proof, we replace $\eZ$ by a semistable model $\eZ'$ of the generic fibre of $\eZ$ dominating $\eZ$. 
\end{proof}

\begin{cor}
\label{t5} Let $\mu$ be any rational number.

i) Every bundle in $\eB_X^\mu$ is semistable of slope $\mu$.

ii) If $E$ is a vector bundle in the category $\eB_X^\mu$ and $F$ is a vector bundle in the category $\eB_X^\nu$, then $E \otimes F$ lies in the category $\eB_X^{\mu + \nu}$. Besides, the internal hom bundle $\mathcal{H}om(E,F)$ lies in the category $\eB_X^{\nu - \mu}$.

iii) The category $\eB_X^\mu$ is abelian. For a point $x \in X (\C_p)$ consider the fibre functor $\omega_x$ on $\eB^0_X$ defined by $\omega_x (E) = E_x$. Then the pair $(\eB^0_X , \omega_x)$ is a neutral Tannakian category over $\C_p$. 

iv) Let $f: X' \rightarrow X$ be a finite morphism of smooth, projective, connected curves over $\oQ_p$ which has degree $\delta$. Then a vector bundle $E$ is in $\eB_X^\mu$ if and only if $f^* E$ is in $\eB^{\delta \mu}_X$. 

v) Any line bundle $\Mh$ of degree $d$ on $X_{\C_p}$ lies in $\eB^d_X$.
\end{cor}

\begin{proof}
i) By \cite{dw1}, theorem 17 and theorem 13, all vector bundles in the category $\eB_X^s$ are semistable of slope zero. Hence our claim follows from theorem \ref{t3}.

ii) By \cite{dw1}, theorem 17 and proposition 9, the category $\eB_X^s$ is closed under tensor products and internal homs. Again the claim follows from theorem \ref{t3}. A more direct proof is also possible using the fact that the class of vector bundles in definition \ref{t1} is closed under tensor products and internal homs. 

iii) By \cite{dw2}, corollary 12.3.4  the category $\eB_X^{ps}$ is abelian. Using theorem \ref{t3}, part iii) we find that for any slope $\mu$ the category 
$\eB_X^\mu$ is an additive full subcategory of the category of all vector bundles on $X_{\C_p}$, which is closed under kernels and cokernels. Hence it is abelian. 

By ii) the category $\eB_X^0$ is also stable under tensor product and internal homs, hence it is a neutral Tannaka category with respect to $\omega_x$. Note that $\omega_x$ is faithful on the category of semistable vector bundles of degree zero on $X_{\C_p}$ and in particular on $\eB^0_X$. Alternatively the faithfulness of $\omega_x$ on $\eB^0_X$ follows as in \cite{dw1} Proposition 30 from the existence of the parallel transport for the bundles in $\eB^0_X$ c.f. Theorem 10 below. 

iv) If $f: X' \rightarrow X$ is a finite morphism of degree $\delta$ as in the claim, then pullback of vector bundles induces a functor $\eB_{X'}^s \rightarrow \eB_X^s$ by \cite{dw1}, theorem 17 and proposition 9. By theorem \ref{t3}, this implies our claim. One can also argue with the functoriality  of the class of vector bundles in definition \ref{t1}. 

v) This follows by taking $\alpha = \id , E = \Mh$ and $L = \Mh^*$ in theorem \ref{t3}, ii).
\end{proof}

Prompted by a question of Faltings in \cite{fa} and encouraged by Corollary \ref{t5} we make the following conjecture:

\begin{conj}
  \label{t6n}
Let $X$ be a smooth projective curve over $\oQ_p$. Every semistable vector bundle of slope $\mu$ on $X_{\C_p}$ lies in $\eB^{\mu}_X$.
\end{conj}

\section{Ramified coverings of curves}

For every variety $Z$ over $\oQ_p$ we denote by $\Pi_1(Z)$ the following topological gruppoid. Its set of objects is $Z(\C_p)$, and for objects $z, z' \in Z(\C_p)$ the set of morphisms $\mbox{Mor}_{\Pi_1(X)}(z,z')$  is the set of isomorphisms of \'etale fibre functors $F_z \rightarrow F_{z'}$. Here $F_z$ is the functor from the category of finite \'etale coverings of $Z$ to the category of finite sets given by $F_z = \Mor_Z (z , \_)$. We call any such morphism an  \'etale path (up to homotopy) from $z$ to $z'$. As a profinite set $\mbox{Mor}_{\Pi_1(X)}(z,z')$ carries a natural topology.
If $z \in Z(\C_p)$ is a base point, we denote by $\pi_1(Z,z)$ the group of closed \'etale paths in $z$, i.e. the algebraic fundamental group of $Z$ with base point $z$. Any morphism of varieties $f: Z \rightarrow Z'$ induces a natural functor $f_\ast: \Pi_1(Z) \rightarrow \Pi_1(Z')$. 

For the purpose of descending the parallel transport we need the following construction:

\begin{construction}\label{t6}
\rm 
  Consider a finite \'etale Galois covering $\alpha : V \to U$ with Galois group $G$ of varieties $U$ and $V$ over $\oQ_p$. Let $B : \Pi_1 (V) \to \Ch$ be a functor to an abelian category $\Ch$. Assume that there is a system of isomorphisms $\varphi_{\sigma} : B \verk \sigma_* \overset{\sim}{\longrightarrow} B$ for $\sigma \in G$ satisfying $\varphi_e = \id$ and $\varphi_{\sigma \tau} = \varphi_{\tau} \verk (\varphi_{\sigma} \tau_*)$. Define a functor $A = A_B : \Pi_1 (U) \to \Ch$ as follows: For $x \in \mathrm{Ob}\, \Pi_1 (U) = U (\C_p)$ set
\[
A (x) = \{ (f_y) \in \prod_{y \in \alpha^{-1} (x)} B (y) \tei \varphi_{\sigma,y} (f_{\sigma y}) = f_y \; \mbox{for all} \; \sigma \in G \; \mbox{and} \; y \in \alpha^{-1} (x) \} \; .
\]
Let $\gamma$ be an \'etale path in $U$ from $x_1 \in U (\C_p)$ to $x_2 \in U (\C_p)$ and let $y_1$ be a point in $V (\C_p)$ with $\alpha (y_1) = x_1$. Then there is a unique path $\delta$ from $y_1$ to some point $y_2 \in V (\C_p)$ over $x_2$ such that $\alpha_* (\delta) = \gamma$. For $\sigma \in G$ the path $\sigma_* (\delta)$ is the unique lift of $\gamma$ starting in $\sigma y_1$ (and ending in $\sigma y_2$). Applying the functor $B$ gives an isomorphism $B (\sigma_* (\delta)) : B (\sigma y_1) \to B (\sigma y_2)$. The product of all $B (\sigma_* (\delta))$ induces an isomorphism
\[
A (\gamma) : \prod_{y \in \alpha^{-1} (x_1)} B (y) = \prod_{\sigma \in G} B (\sigma y_1) \longrightarrow \prod_{\sigma \in G} B (\sigma y_2) = \prod_{y \in \alpha^{-1} (x_2)} B (y) \; .
\]
It maps $A (x_1)$ to $A (x_2)$ by the naturality of the $\varphi_{\sigma}$'s and it does not depend on the choice of the point $y_1$ over $x_1$. Explicitely $A (\gamma)$ is given as follows. For an element $(f_{\sigma y_1})_{\sigma \in G}$ of $A (x_1)$ we have
\begin{equation}
  \label{eq:1}
  A (\gamma) ((f_{\sigma y_1})_{\sigma \in G}) = (B (\sigma_* (\delta)) (f_{\sigma y_1}))_{\sigma \in G} \quad \mbox{in} \; A (x_2) \; .
\end{equation}
Here $B (\sigma_* (\delta)) (f_{\sigma y_1})$ is in $B (\sigma y_2)$.

Having thus constructed a functor $A : \Pi_1 (U) \to \Ch$ we define a canonical isomorphism of functors:
\begin{equation}
  \label{eq:2}
  \Phi : A \verk \alpha_* \longrightarrow B \; .
\end{equation}
Namely, for every $z \in V (\C_p)$ with $\alpha (z) = x$ the projection map 
\[
\prod_{y \in \alpha^{-1} (x)} B (y) \to B (z)
\]
induces an isomorphism $\Phi (z) : A (x) \to B (x)$. This follows from the cocycle condition for the $\varphi_{\sigma}$'s above. For all $\sigma \in G$ the following diagram commutes:
\begin{equation}
  \label{eq:3}
  \xymatrix{A \verk \alpha_* \verk \sigma_* \ar[r]^{\Phi \verk \sigma_*} \ar[d]^{\wr \, \mathrm{can}} & B \verk \sigma_* \ar[d]^{\varphi_{\sigma}} \\
A \verk \alpha_* \ar[r]^{\Phi} & B \; .}
\end{equation}
If $\Ch$ is a topological abelian category and $B$ above a continuous functor, then the functor $A$ is continuous as well.
\end{construction}

\begin{rem}
  The explicit construction of the functor $A = A_B$ and our later calculations with it are somewhat clumsy. However since we have to compute group actions  later and since non-trivial actions could hide in isomorphisms we felt more secure with this pedestrian setup. 
\end{rem}

As before we fix a smooth, projective and connected curve $X$ over $\oQ_p$ and write $X_{\C_p}$ for the base change of $X$ to $\C_p$. 
We consider a finite Galois covering
\[\alpha: Y \longrightarrow X\]
of $X$ with Galois group $G$, by a smooth, projective and connected curve $Y$ over $\oQ_p$. 
By $S \subset X(\oQ_p)$ we denote the ramification locus of $\alpha$, and we put $T = \alpha^{-1}(S) \subset Y(\oQ_p)$. 
Moreover, let  $U = X \backslash S$ and $V = Y \backslash T$. We denote by $j: V \hookrightarrow Y$ the corresponding open immersion.

Note that the \'etale Galois covering $\alpha|_V: V \rightarrow U$ gives rise to an exact sequence
\begin{equation} \label{eq:4} 
1 \rightarrow \pi_1(V,y_0) \rightarrow \pi_1(U, x_0) \rightarrow G \rightarrow 1,
\end{equation}
where $y_0 \in V(\C_p)$ and $x_0 = \alpha(y_0)$ are base points.

Throughout this section we fix a vector bundle $F$ on $Y_{\C_p}$ carrying a $G$-action, i.e. we have isomorphisms
\[a_\sigma: \sigma^\ast F \iso F \mbox{ for all }\sigma \in G\]
satisfying $a_{\sigma \tau} = a_\tau \circ \tau^\ast(a_\sigma)$ and $a_{e} = \mbox{id}$. 
Assume that $F$ lies in the category $\eB_Y^{ps}$. 

In \cite{dw1}, Theorem 36 and remark, we defined a continuous functor $\rho_F$ from $\Pi_1(Y)$ to the category $\mbox{\bf Vec}_{\C_p}$ of finite dimensional $\C_p$-vector spaces. For every $z \in Y(\C_p)$ we have $\rho_F(z) = F_z$, where $F_z$ is the fibre of $F$ in $z$. For a morphism $\gamma$ from $y_1$ to $y_2$, i.e. an \'etale path from $y_1$ to $y_2$, the isomorphism $\rho_{F}(\gamma) : F_{y_1} \rightarrow F_{y_2}$  is called the parallel transport along $\gamma$. The association $F \mapsto \rho_F$ is functorial. In particular, if $\sigma$ is an element of $G$, there is a commutative diagram 
\begin{equation} \label{eq:5}
\xymatrix{
F_{\sigma( y_1)}  \ar[r]^{a_{\sigma, y_1}} \ar[d]_{\rho_{\sigma^\ast F}(\gamma)} & F_{y_1} \ar[d]^{\rho_F(\gamma)} \\
F_{\sigma( y_2)} \ar[r]_{a_{\sigma, y_2}}& F_{y_2} .
}
\end{equation}
Since our construction is compatible with pullbacks by \cite{dw1}, theorem 36, we have $ \rho_F \circ \sigma_\ast = \rho_{\sigma^\ast F}$, where $\sigma_\ast: \Pi_1(Y) \rightarrow \Pi_1(Y)$ is induced by $\sigma$. The system of isomorphisms of continuous functors from $\Pi_1 (V)$ to $\uVec_{\C_p}$
\[
\varphi_{\sigma} : (\rho_F \verk j_*) \verk \sigma_* = (\rho_F \verk \sigma_*) \verk j_* = \rho_{\sigma^* F} \verk j_* \overset{\overset{a_{\sigma} \verk j_*}{\sim}}{\longrightarrow} \rho_F \verk j_*
\]
satisfies the cocycle condition of construction \ref{t6}. Hence we get a continuous functor
\begin{equation}
  \label{eq:6}
  \rho = A_{\rho_F \verk j_*} : \Pi_1 (U) \longrightarrow \uVec_{\C_p}
\end{equation}
and a canonical isomorphism $\Phi : \rho \verk \alpha_* \overset{\sim}{\longrightarrow} \rho_F \verk j_*$ satisfying $\varphi_{\sigma} \verk (\Phi \verk \sigma_*) = \Phi$ for every $\sigma \in G$.

In particular, for every base point $x_0 \in U (\C_p)$ we get a continuous representation 
\begin{equation} \label{eq:7}
\rho: \pi_1(U,x_0) \rightarrow \mathrm{GL}(\rho(x_0)).
\end{equation}

Now we want to show that $\rho$ factors over a certain quotient of $\pi_1(U, x_0)$ by
studying the monodromy along \'etale paths around points in the ramification divisor $S$. 
For every closed point $y$ in $Y$ we denote by $G_y$ the subgroup of $G$ consisting of all elements fixing $y$. For $\sigma \in G_y$ the map $a_{\sigma} : \sigma^* F \overset{\sim}{\longrightarrow} F$ has fibre $a_{\sigma , y} : F_y = (\sigma^* F)_y \to F_y$, an automorphism of $F_y$. The map $G_y \to \GL (F_y)$ sending $\sigma$ to $a_{\sigma , y}$ is an antihomomorphism. 

\begin{lemma}
\label{t7}
Let $t$ be a point in $T$ (viewed in $T (\C_p)$) and assume that for all $\sigma \in G_t$ the element $a_{\sigma,t} \in \mathrm{GL}(F_t)$ is central, i.e. there exists a character
\[\chi_t: G_t \rightarrow \mathbb{C}_p^\times\] satisfying
$a_{\sigma, t} = \chi_t(\sigma) \mathrm{id}_{F_t}$. We fix a base point $x_0$ in $U(\C_p)$ and a preimage $y_0$ of $x_0$ in $V(\C_p)$. Let $\gamma_0$ be an element in $\pi_1(U, x_0)$ mapping via \eqref{eq:4} to an element $\tau \in G_t$, and let $\gamma$ be an \'etale path in $V$ with starting point $y_0$ that lifts $\gamma_0$. Its endpoint is therefore $\tau (y_0)$. 
If there exists an \'etale path $\delta$ in $Y$ from $y_0$ to $t$ satisfying
\[(\ast) \quad \quad \rho_{F}(\tau_\ast (\delta) j_\ast (\gamma) \delta\inv) = \mathrm{id}_{F_t},\]
then we have
\[\rho(\gamma_0) = \chi_t(\tau) \mathrm{id}_{\rho(x_0)}.\]
If $(\ast)$ holds for one \'etale path $\delta$ from $y_0$ to $t$, it holds for all such paths.
\end{lemma}

\begin{rem}
  We will verify condition $(\ast)$ in the proof of theorem \ref{t9} below using Grothendieck's comparison theorem between algebraic and topological fundamental groups. We do not know an algebraic argument showing that $(\ast)$ is satisfied.
\end{rem}

\begin{proof}
Since $\rho_F$ is a functor, we find
\begin{eqnarray*}
\mathrm{id_{F_t}} & = & \rho_F(\tau_\ast(\delta) j_\ast (\gamma) \delta\inv )  \\
& = & \rho_F(\tau_\ast (\delta)) \circ \rho_F(j_\ast (\gamma)) \circ \rho_F(\delta\inv),
\end{eqnarray*}
hence $\rho_F(j_\ast (\gamma)) = \rho_F(\tau_\ast (\delta))\inv \circ \rho_F(\delta)$.
Now according to \eqref{eq:5} we have
\[
\rho_F(\tau_\ast (\delta)) = \rho_{\tau^* F} (\delta) = 
a_{\tau, t}\inv \circ \rho_F( \delta) \circ a_{\tau, y_0}.
\]
By assumption, we have $a_{\tau, t} = \chi_t(\tau) \mathrm{\id}_{F_t}$, hence
\[\rho_F(j_\ast (\gamma))= \rho_F(\tau_\ast (\delta))\inv \circ \rho_F(\delta) = a_{\tau,y_0}\inv \circ  \rho_F(\delta)\inv \circ a_{\tau,t}\circ \rho_F(\delta) = \chi_t(\tau) a_{\tau, y_0}\inv.\]
Hence for every element $(f_y)$ of 
\[
\rho(x_0)= \{(f_y) \in \prod_{y \in \alpha^{-1} (x_0)} F_{y} \tei a_{\sigma,y} (f_{\sigma y}) = f_y \mbox{ for all }\sigma \in G\}
\]
we find 
\[\rho_F(j_\ast (\gamma))(f_{y_0}) = \chi_t(\tau) a_{\tau,y_0}\inv (f_{y_0}) = \chi_t (\tau) f_{\tau y_0}.\]
This implies that $\rho(\gamma_0) = \chi_t(\tau) \mathrm{id}_{\rho(x)}$. Namely, according to \eqref{eq:1} the component of $\rho (\gamma_0) (f_{\sigma y_0})_{\sigma \in G}$ in $F_{\sigma \tau y_0}$ is given by
\begin{eqnarray*}
  \rho_F (\sigma_* (j_* (\gamma))) (f_{\sigma y_0}) & = & a^{-1}_{\sigma , \tau y_0} (\rho_F (j_* (\gamma)) (f_{y_0})) \\
& = & a^{-1}_{\sigma , \tau y_0} (\chi_t (\tau) f_{\tau y_0}) \\
& = & \chi_t (\tau) a^{-1}_{\sigma, \tau y_0} (f_{\tau y_0}) \\
& = & \chi_t (\tau) f_{\sigma \tau y_0} \; .
\end{eqnarray*}
This shows that as claimed we have
\[
\rho (\gamma_0) ((f_y)) = \chi_t (\tau) (f_y)
\]
where $(f_y) = (f_y)_{y \in \alpha^{-1} (x_0)}$.

Finally, we check that the condition $(\ast)$ is independent of the choice of $\delta$. If $\delta'$ is a second \'etale path in $Y$ from $y_0$ to $t$, then $\eta = \delta' \delta\inv$ is an element of $\pi_1(Y,t)$. Since $a_{\tau,t}$ is central, we deduce from \eqref{eq:5} that we have
\[\rho_F(\tau_\ast(\eta)) = \rho_{\tau^* F} (\eta) = a_{\tau, t}\inv \circ \rho_F(\eta) \circ a_{\tau, t} = \rho_F(\eta).\]
Hence, using condition $(\ast)$ for $\delta$, we find
\begin{eqnarray*}
 \rho_F(\tau_\ast(\delta') j_\ast (\gamma) (\delta')\inv) & = 
& \rho_F(\tau_\ast(\eta \delta) j_\ast (\gamma) \delta\inv \eta\inv) \\
~ & = & \rho_F(\tau_\ast(\eta)) \rho_F(\eta)\inv\\
~ & = & \mathrm{id}_{F_t},
\end{eqnarray*}
hence $(\ast)$ also holds for $\delta'$. 
\end{proof}

Now we use the transcendental theory of the fundamental group. Fix an embedding $\C_p \subset \C$ of abstract fields and set $X_{\C} = X \otimes_{\oQ_p} \C$ etc. By $X^{\ann} = (X_{\C})^{\ann}$ etc. we denote the associated Riemann surfaces. Any point $x_0 \in U (\C_p)$ then defines points $x_0$ in $U_{\C} (\C)$ and $U^{\ann}$. By \cite{sga1} expos\'e X, Corollaire 1.8 and expos\'e XII, Corollaire 5.2 we have canonical isomorphisms
\[
\hat{\pi}_1 (U^{\ann}, x_0) \overset{\sim}{\longrightarrow} \pi_1 (U_{\C} , x_0) \overset{\sim}{\longrightarrow} \pi_1 (U , x_0) \; .
\]
Here $\hat{\pi}_1 (U^{\ann} , x_0)$ denotes the profinite completion of the topological fundamental group $\pi_1 (U^{\ann} , x_0)$. Let $g$ be the genus of $X$. For every $s \in S \cong S_{\C}$ we choose a closed loop $\gamma_s$ with basepoint $x_0$ in $U^{\ann}$ such that $\gamma_s$ is contained in a contractible open subset of $U^{\ann} \cup \{ s \}$ and winds once around $s$, counterclockwise with respect to the orientation. Moreover, for all $i = 1 , \ldots , g$ let $\gamma_i$ and $\gamma_{g+i}$ be closed loops in $U^{\ann}$ with basepoint $x_0$ such that $\pi_1 (X^{\ann} , x_0)$ is the group generated by $\gamma_1 , \ldots , \gamma_{2g}$ subject to the single relation
\[
\prod^g_{i=1} [\gamma_i , \gamma_{g+i}] = 1 \; .
\]
Then $\pi_1 (U^{\ann} , x_0)$ is the group generated by $\gamma_1 , \ldots , \gamma_{2g}$ and all $\gamma_s$ for $s \in S$ subject to the relation (where $S$ carries a suitable ordering)
\[
\prod_{s \in S} \gamma_s \prod^g_{i=1} [\gamma_i , \gamma_{g+i}] = 1 \; .
\]
By abuse of notation we denote by $\gamma_i$ and $\gamma_s$ also the corresponding elements in the algebraic fundamental group under the map
\[
\pi_1 (U^{\ann} , x_0) \longrightarrow \hat{\pi}_1 (U^{\ann} , x_0) \cong \pi_1 (U , x_0) 
\]
and similarly for $\pi_1 (X , x_0)$. Then $\pi_1 (X , x_0)$ is the profinite group generated by all $\gamma_i$ subject to the relation $\prod^g_{i=1} [\gamma_i , \gamma_{i+g}] = 1$. Set
\[
N = \Ker (\pi_1  (U , x_0) \longrightarrow \pi_1 (X , x_0)) \; .
\]
Then $N$ is the closed normal subgroup of $\pi_1 (U  , x_0)$ generated by the elements $\gamma_s$ for $s \in S$. Note that $\gamma_s$ generates an inertia group $I_s$ at $s$. Let $[N , \pi_1 (U , x_0)]$ be the closed normal subgroup of $\pi_1 (U , x_0)$ generated by all commutators $[n, \gamma]$ for $\gamma \in \pi_1 (U , x_0)$ and $n \in N$ or equivalently for $n = \gamma_s$ for $s \in S$. Setting
\[
\tilde{N} = N / [N , \pi_1 (U , x_0)] \quad \mbox{and} \quad  \tilde{\Gamma} = \pi_1 (U , x_0) / [N , \pi_1 (U , x_0)] \; ,
\]
we get a canonical central extension
\[
1 \to \tilde{N} \to \tilde{\Gamma} \to \pi_1 (X , x_0) \to 1 \; .
\]
All inertia groups at $s \in S$ in $N$ have the same image $\tilde{I}_s$ in $\tN$ because $\tN$ is central in $\tilde{\Gamma}$. Moreover, by Kummer theory each $\tilde{I}_s$ is canonically isomorphic to $\hZ (1) = \varprojlim_n \mu_n$ where $\mu_n = \mu_n (\C_p)$. The image $\tilde{\gamma}_s$ of $\gamma_s$ in $\tilde{N}$ is a generator of $\tilde{I}_s$. For $g \ge 1$ we have
\[
\tN = \prod_{s \in S} \tilde{I}_s = \hZ (1)^{|S|} \; .
\]
For $g = 0$ the group $\tN$ is the quotient of $\prod_{s \in S} \tilde{I}_s$ by the single relation $\prod_{s \in S} \tilde{\gamma}_s = 1$ and hence it is isomorphic to $\hZ (1)^{|S|-1}$. \\
Now fix integers $d_s \ge 1$ for $s \in S$ and let $\Gamma$ be the quotient of $\tilde{\Gamma}$ by the relations $\gamma^{d_s} = 1$ for all $\gamma \in \tilde{I}_s$, ($\gamma = \gamma_s$ suffice) and $s \in S$.

The image $\overline{I}_s$ of $\tilde{I}_s$ in $\Gamma$ is canonically isomorphic to $\mu_{d_s}$ and generated by $\overline{\gamma}_s$, the image of $\tilde{\gamma}_s$ in $\overline{\Gamma}$.  There is an exact sequence
\begin{equation}
  \label{eq:8}
  1 \to \oN \to \Gamma \to \pi_1 (X , x_0) \to 1
\end{equation}
where 
\[
\oN = \prod_{s \in S} \overline{I}_s = \prod_{s \in S} \mu_{d_s} \quad \mbox{if} \; g \ge 1 \; .
\]
For $g = 0$ the group $\Gamma$ is the quotient of $\prod\limits_{s \in S} \overline{I}_s$ by the single relation $\prod\limits_{s \in S} \overline{\gamma}_s = 1$.


We now have the following result:
\begin{theorem}
  \label{t9}
Consider as above a Galois covering $\alpha : Y \to X$ with group $G$, ramification loci $S \subset X$ and $T = \alpha^{-1} (S) \subset Y$ with complements $U = X \ohne S$ and $V = Y \ohne T$. For a vector bundle $F$ in $\eB^{ps}_Y$ with $G$-action given by $a_{\sigma} : \sigma^* F \to F$ for every $\sigma \in G$ and for a point $x_0 \in U (\C_p)$ consider the representation \eqref{eq:7} obtained by descending $\rho_F$: 
\[
\rho : \pi_1 (U , x_0) \longrightarrow \GL (\rho (x_0)) \; .
\]
Assume that for every point $t \in T$ and every $\sigma \in G_t$ the element $a_{\sigma , t} \in \GL (F_t)$ is central. Let $\chi_t : G_t \to \mu \subset \C^*_p$ be the corresponding character satisfying $a_{\sigma , t} = \chi_t (\sigma) \id_{F_t}$. Given $s \in S$, let $t \in T$ be a point with $\alpha (t) = s$. Then there exists a generator $\tau$ of $G_t$ such that $\rho (\gamma_s)$ is the automorphism 
\[
\rho (\gamma_s) = \chi_t (\tau) \id : \rho (x_0) \longrightarrow \rho (x_0) \; .
\]
The order $o_s$ of the character $\chi_t$ depends only on $s$ and divides the ramification index $e_s$ of the point $s$. Fix integer multiples $d_s \ge 1$ of $o_s$ for $s \in S$ and let $\Gamma$ be the quotient of $\pi_1 (U, x_0)$ constructed above using the numbers $d_s$ for $s \in S$. Then the representation $\rho$ factors over $\Gamma$ and induces a continuous representation
\[
\rho : \Gamma \longrightarrow \GL (\rho (x_0)) \; .
\]
\end{theorem}

\begin{proof}
We look at the analytic covering $\alpha^{an}: Y^{an} \rightarrow X^{an}$ induced by $\alpha$. It is unramified over $U^{an}$. For any $s \in S$ choose a contractible open non-compact neighbourhood $V_s$ of $\gamma_s$ in $U^{an} \cup \{s\}$ which contains $s$. Shrinking $V_s$ if necessary we may assume that $V_s$ is not biholomorphic to $\C$. By Riemann's mapping theorem there is then a biholomorphic map from $V_s$ to the open unit disc $E$, mapping $s$ to $0$. 
Let $W$ be the connected component of $(\alpha^{an})\inv V_s$ containing $t$. 
By \cite{fo}, Satz 5.11, there exists a biholomorphic map $\varphi: W \rightarrow E$ such that the diagram
\[
\xymatrix{
W \ar[r]^{\overset{\varphi}{\sim}} \ar[d]_{\alpha^{an}} &  E \ar[d]^{z \mapsto z^{e_s}}\\
V_s \ar[r]^{\sim} & E
}
\]
commutes.
In particular we have $\varphi (t) = 0$. Choose a point $y_0 \in W$ mapping to $x_0$ and lift $\gamma_s$ to a path $\gamma$ in $W \backslash \{t\}$ with starting point $y_0$. Its endpoint is a point in $W$  of the form $\tau y_0$ for some $\tau \in G$ of order $e_s$. Note here that $\gamma_s$ winds around $s$ once. Hence $\tau$ maps the connected component $W$ to itself, which implies $\tau t = t$, i.e. $\tau$ is contained in $G_t$. Since $|G_t| = e_s$ it follows that $\tau$ generates $G_t$. 

Choose a path $\delta$ in $W$ from $y_0$ to $t$. Then $\tau (\delta)  \gamma \delta\inv$ is a closed path in $W$. Since $W$ is biholomorphic to the unit disc $E$, this path is null-homotopic in $W$, and hence also in $Y^{an}$. The topological path $\delta$ induces an \'etale path in $Y$ connecting $y_0$ and $t$. By abuse of notation, we denote it also by $\delta$. It follows that the composition $\tau_* (\delta) j_* (\gamma) \delta^{-1}$ is trivial in the algebraic fundamental group $\pi_1(Y, t)$. Hence we can apply lemma \ref{t7} to conclude that 
\[\rho(\gamma_s) = \chi_t(\tau) \mathrm{id}_{ \rho(x_0)}.\]
In particular, $\rho(\gamma_s)$ is central in $\mathrm{GL}(\rho(x_0))$ for every $s \in S$. Therefore \\
$\rho: \pi_1(U, x_0) \rightarrow \mathrm{GL}(\rho(x_0))$ factors over $\tilde{\Gamma}$. Given two points $t , t'$ over $s$ there is some $\pi \in G$ with $\pi t = t'$. The cocycle relation for the $a_{\sigma}$'s implies that we have a commutative diagram
\[
\xymatrix{
G_t \ar[rr]^{\sim} \ar[dr]_{\chi_t} && G_{t'}\ar[dl]^{\chi_{t'}} \\
 & \mu &
}
\]
where the horizontal map is given by conjugation with $\pi$. Hence the orders of $\chi_t$ and $\chi_{t'}$ are the same and divide $e_s = |G_t|$. We have $\rho (\gamma_s^{d_s})= \rho (\gamma_s)^{d_s} = \chi_t(\tau)^{d_s} = 1$ since $o_s \tei d_s$, so that indeed $\rho$ induces a representation
\[\rho: \Gamma \longrightarrow \mathrm{GL}(\rho(x_0)).\]
Note here that $\tilde{\gamma}_s$ generates $\tilde{I}_s$ for every $s \in S$.
\end{proof}

Theorem \ref{t9} will be used in the next section to define representations for certain bundles of non-zero slope. For now we use it to define a parallel transport for bundles in $\eB_X^0$, the category of vector bundles of slope zero on $X$ with potentially strongly semistable reduction in the sense of definition \ref{t2}.
\begin{theorem}
\label{t10}
Let $E$ be a vector bundle on $X_{\C_p}$ contained in the category $\eB_X^0$. The construction below gives a continuous functor
\[\rho_E: \Pi_1(X) \longrightarrow \mathrm{Vec}_{\C_p}\]
from the \'etale fundamental groupoid $\Pi_1(X)$ to the category of finite-dimensional $\C_p$-vector spaces,
mapping every $x \in X(\C_p)$ to the fibre $E_x$. In other words, there is a functorial parallel transport on $E$ along \'etale paths on $X$. For bundles in $\eB^0_X$ parallel transport is compatible with $\otimes, \uHom , \Gal (\oQ_p / \Q_p)$-conjugation and pullbacks via morphisms of smooth projective curves over $\oQ_p$. In particular, for all $x \in X (\C_p)$, we obtain a representation $\rho_{E,x}: \pi_1(X,x) \rightarrow \mathrm{GL}(E_x)$ of the \'etale fundamental group satisfying the corresponding compatibilities. For bundles in $\eB^{ps}_X$ the parallel transport is the same as the one defined in \cite{dw1}. 
\end{theorem}

\begin{proof} 
Since $E$ lies in $\eB_X^0$, by theorem \ref{t3} there exists a finite morphism \\
$\alpha: Y \rightarrow X$ of smooth, projective, connected curves over $\oQ_p$, such that $F = \alpha_{\C_p}^\ast E$ lies in $\eB_Y^{ps}$. We may assume that $\alpha$ is Galois with group $G$. Then we have canonical identifications $a_{\sigma} : \sigma^* F = F$ satisfying the cocycle conditions and hence we get the descended functor \eqref{eq:6}:
\[
\rho = A_{\rho_F \verk j_*} : \Pi_1 (U) \longrightarrow \uVec_{\C_p} \; .
\]
For $x \in \Ob \; \Pi_1 (U) = U (\C_p)$ we have
\begin{eqnarray*}
  \rho (x) & = & \{ (f_y) \in \prod_{\alpha (y) = x} F_y \tei a_{\sigma , y} (f_{\sigma y}) = f_y \; \mbox{for all} \; \sigma \in G \} \\
& = & \{ (f_y) \in \prod_{\alpha (y) = x} E_x \tei f_{\sigma y} = f_y \; \mbox{for all} \; \sigma \in G \} \; .
\end{eqnarray*}
Hence we may identify $\rho (x)$ and $E_x$. Clearly, for $\sigma \in G_t$ 
 the action of $a_{\sigma, t}$ on $F_t$ is given by the trivial character. For any $x_0 \in U (\C_p)$ applying theorem \ref{t9}, we find that $\rho(\gamma_s) = 1$. Therefore the representation $\pi_1(U,x_0) \longrightarrow \mathrm{GL}(E_{x_0})$ induced by $\rho$ factors through $\pi_1(X, x_0)$. 
 
Let $\lambda : U \hookrightarrow X$ be the open immersion. Recall that $j: V \hookrightarrow Y$ denotes the corresponding map on the covering.
Now let $\gamma$ be any \'etale path on $X$ from $x_1\in X (\C_p)$ to $x_2 \in X (\C_p)$. We claim that $\gamma$ can be lifted to an \'etale path in $Y$. Choose preimages $y_1$ and $y_2$ of $x_1$ and $x_2$ in $Y(\C_p)$. By \cite{sga1},expos\'e V, no. 7,  there exists an \'etale path $\epsilon$ in $Y$ connecting $y_2$ and $y_1$. Then $\alpha_\ast(\epsilon) \gamma$ is a closed \'etale path in $X$ based in $x_1$. If $x_1$ lies in $U(\C_p)$, then $\alpha_\ast(\epsilon) \gamma = \lambda _\ast(\delta)$ for some closed \'etale path $\delta$ on $U$, since $\pi_1(U, x_1) \rightarrow \pi_1(X,x_1)$ is surjective. Since $\alpha: V \rightarrow U$ is \'etale, there exists an \'etale path $\delta'$ on $V$ lifting $\delta$ with endpoint $y_1$. Then 
\[
\alpha_\ast(\epsilon\inv j_\ast(\delta')) = \alpha_\ast (\epsilon)\inv  \lambda_\ast \alpha_\ast (\delta') = \alpha_\ast (\epsilon)\inv \lambda_\ast (\delta)  = \alpha_\ast (\epsilon)\inv  \alpha_\ast(\epsilon) \gamma = \gamma.
\]
Hence any path $\gamma$ starting in a point of $U (\C_p)$ can be lifted to an \'etale path in $Y$. 
If $x_1$ does not lie in $U(\C_p)$, we choose a point $x_0 \in U(\C_p)$, a preimage $y_0$ of $x_0$ in $V(\C_p)$ and an \'etale path $\delta_1$ from $y_0$ to $y_1$ in $Y$. Then $\gamma \alpha_\ast (\delta_1)$ is an \'etale path in $X$ starting in $x_0$. As we have just shown, it can be lifted to an \'etale path in $Y$. This implies that  $\gamma$ can be lifted to an \'etale path in $Y$, too.

We claim that there is a representation
\[\rho_E: \Pi_1(X) \longrightarrow \mathrm{Vec}_{\C_p}\]
such that $\rho$ above factors as
\[
\rho : \Pi_1 (U) \to \Pi_1 (X) \xrightarrow{\rho_E} \uVec_{\C_p} \; .
\]
We define $\rho_E$ by setting $\rho_E(x) = E_x$ for all $x \in X(\C_p)$ and by 
\[\rho_E(\gamma) = \rho_F(\delta),\]
if $\delta$ is an arbitrary \'etale path in $Y$ satisfying $\alpha_\ast (\delta) = \gamma$. Here we identify $F_y$ with $E_x$ for any point $y$ above $x$.

We have already shown that such a lift $\delta$ always exists. It remains to check that the definition  is independent of the choice of $\delta$. Let $\gamma$ be an \'etale path in $X$ starting in $x_1 \in U(\C_p)$, and assume that $\delta_1$ and $\delta_2$ are \'etale paths in $Y$ lifting $\gamma$. We denote by $y_1 \in V(\C_p)$ the  starting point of $\delta_1$. Replacing $\delta_2$ by $\sigma_\ast \delta_2$ for a suitable $\sigma \in G$, we may assume that $\delta_1$ and $\delta_2$ have the same endpoint. Hence $\delta_2^{-1} \delta_1$ is an \'etale path in $Y$ connecting two points in $V$. Since for all $y_0 \in V(\C_p)$ the map $\pi_1(V,y_0) \rightarrow \pi_1(Y, y_0)$ is surjective, there exists an \'etale path $\epsilon$ in $V$ such that $j_\ast \epsilon = \delta_2\inv \delta_1$. Hence
\[
\rho_F (\delta_2)^{-1} \rho_F (\delta_1) = \rho_F(\delta_2\inv \delta_1) = \rho_F(j_\ast (\epsilon)) = \rho(\alpha_\ast (\epsilon))
\]
by definition of $\rho$. Now 
$\alpha_\ast (\epsilon)$ is a closed \'etale path in $U$ satisfying $\lambda_\ast \alpha_\ast (\epsilon) = \alpha_\ast j_\ast (\epsilon) = \alpha_\ast (\delta_2\inv \delta_1) = 1$. Therefore it lies in the kernel of $\pi_1(U, x_1) \rightarrow \pi_1(X, x_1)$, which implies $\rho(\alpha_\ast (\epsilon)) = 1$, as we have shown at the beginning of the proof. 

It remains to treat the case that the starting point $x_1$ of $\gamma$ does not belong to $U$. Assume that $\delta_1$ and $\delta_2$ are \'etale paths in $Y$ lifting $\gamma$. Let $y_1 \in Y(\C_p)$ be the starting point of $\delta_1$. Replacing $\delta_2$ by $\sigma_\ast (\delta_2)$ for a suitable $\sigma \in G$, we can assume that $\delta_2$ also starts in $y_1$. Choose a point $x_0 \in U(\C_p)$, a preimage $y_0 \in V(\C_p)$ of $x_0$ and an \'etale path $\epsilon$ in $Y$ from $y_0$ to $y_1$. Then $\delta_1 \epsilon$ and $\delta_2 \epsilon$ are paths in $Y$ lifting the \'etale path $\gamma \alpha_\ast \epsilon$. Since $\gamma \alpha_\ast \epsilon$ starts in $x_0 \in U(\C_p)$, we have already shown that $\rho_F(\delta_1 \epsilon) = \rho_F(\delta_2 \epsilon)$, which implies $\rho_F(\delta_1)= \rho_F(\delta_2)$.

Using the corresponding results from \cite{dw1} it follows that $\rho_E$ is a functor which is functorial in $E$ and satisfies the stated compatibilities.
\end{proof}

\begin{remark}
  \label{t10n}
In the construction of the functor $\rho = A_{\rho_F \circ j_*}$ in \eqref{eq:6} we had assumed that $F$ was a vector bundle in $\eB^{ps}_Y$ in order to have the parallel transport $\rho_F$. Using theorem \ref{t10} one can consider the functor $\rho = A_{\rho_F \circ j_*}$ also for the bundles in the a priori bigger category $\eB^0_Y$ and lemma \ref{t7} and theorem \ref{t9} continue to hold with identical proofs, if $F$ is in $\eB^0_Y$.
\end{remark}

\section{Vector bundles of non-vanishing slope}
As usual, let $X$ be a smooth, projective and connected curve over $\oQ_p$.
We have seen in the last section that bundles in $\eB_X^0$ admit functorial isomorphisms of parallel transport. 
In particular, they give rise to $p$-adic representations of the fundamental group $\pi_1(X, x_0)$. 
Now we will look at bundles in $\eB_X^\mu$ for arbitrary slopes $\mu = d/r$. In general we cannot expect $p$-adic representations of the fundamental group. However, similarly to the classical case we will define for any bundle $E$ in $\eB_X^\mu$ a representation $\rho_{E,x_0}$  of a certain central extension $\Gamma_r$ of $\pi_1(X, x_0)$, c.f. \cite{nase}, proposition 6.2 and \cite{ab}, section 6. We hope to study these representations elsewhere.

In \cite{nase} the construction uses an infinite extension of the base Riemann surface which is ramified in exactly one point. In our algebraic situation we might naively try to use $\Z / r$-coverings which are ramified in one point of $X$ only. However such coverings do not exist for $r \ge 2$. Coverings with group $\Z / r$ which are ramified in exactly two points do exist though and we can use them for our purposes.

\begin{construction}
\label{t11}
Fix an integer $r \geq 1$ and two different points $s, s_0$ in $X(\oQ_p)$. Below we construct finite \'etale coverings $\alpha: Y \longrightarrow X$ by smooth, projective, connected $\oQ_p$-curves $Y$ such that 
\begin{compactitem}[$\bullet$]
\item $\alpha$ is Galois with Galois group $G = \mu_r$,
\item  $\alpha$  is unramified outside $\{{s},{s_0}\}$ and 
\item $\alpha$ is totally ramified over ${s}$ and ${s_0}$. 
\end{compactitem}
\end{construction}

\begin{proof}
Put $U = X \backslash \{{s}, {s_0}\}$ and consider the following relative exact sequence of \'etale cohomology groups together with the trace isomorphism:
\[
H^1(X, \mu_r) \longrightarrow H^1(U, \mu_r) \stackrel{\delta}{\longrightarrow} H^2_{\{{s}, {s_0}\}}(X, \mu_r) \longrightarrow H^2(X, \mu_r) \overset{\mathrm{tr}}{\cong} \Z / r.\]
Besides we have: 
\[
H^2_{\{{s}, {s_0} \}}(X, \mu_r) = H^2_{{s}}(X, \mu_r) \oplus H^2_{{s_0}}(X, \mu_r) \simeq \Z/ r \oplus \Z/r \; ,
\]
and these identifications make the following diagram commutative
\[
\xymatrix{
H^1(U, \mu_r) \ar[r]^-{\delta}& H^2_{\{{s},{s_0}\}}(X,\mu_r) \ar[r] \ar[d]_{\simeq} &  H^2(X, \mu_r) \ar[d]^{\simeq}\\
& \Z / r \oplus \Z / r \ar[r]^-{\Sigma} & \Z / r
}
\]
Fix a point $x_0 \in U (\C_p)$ and choose a class $c \in H^1(U, \mu_r) = \Hom (\pi_1 (U , x_0) , \mu_r)$ such that $\delta(c)$ corresponds to $(1,-1)\in \Z/r \Z \oplus \Z / r \Z$. The image of the homomorphism
\[ c : \pi_1(U, x_0) \longrightarrow \mu_r \]
is a subgroup $\mu_t$ of $\mu_r$ for some $t$ dividing $r$. Then $t c$ is trivial, which implies $(t, -t) = t \delta(c) = 0$ in $\Z / r \Z \times \Z / r \Z$. Hence $t = r$, i.e. $c$ is surjective. Its kernel defines a connected Galois \'etale covering $\alpha: V \rightarrow U$ with Galois group $\mu_r$. Let $Y$ be the normalization of $X$ in the function field of $V$. Then $\alpha$ extends uniquely to a finite morphism $\alpha: Y \rightarrow X,$ where $Y$ is a connected, smooth and projective curve over $\oQ_p$. 

Under pullback by $\alpha$, the class $c$ vanishes in $H^1(V, \mu_r)$. Hence $\delta(\alpha^\ast (c))$ is zero in 
\[
H^2_{\alpha^{-1}{\{s, s_0\}}}(Y, \mu_r) = H^2_{\alpha^{-1} (s)} (Y , \mu_r) \oplus H^2_{\alpha^{-1} (s_0)} (Y , \mu_r).
\] 
Denoting by $e$ the common ramification index of the points in $\alpha^{-1} (s)$ we have a commutative diagram with $\Delta_e (a) = (ea , \ldots , ea)$:
\[
\xymatrix{
H^2_{\alpha^{-1} (s)} (Y , \mu_r) \ar@{=}[r] & \bigoplus\limits_{t \in \alpha^{-1} (s)} H^2_t (Y , \mu_r) \ar@{=}[r] & \bigoplus\limits_{t \in \alpha^{-1} (s)} \Z/r \\
H^2_s (X , \mu_r) \ar[u]_{\alpha^*} \ar[rr]^{\sim} && \Z / r \ar[u]_{\Delta_e}
}
\]
Since $\delta (\alpha^* (c))$ vanishes, we have $\Delta_e (1) = 0$ and hence $e = r$. Thus $\alpha$ is totally ramified over $s$ and by the same argument also over $s_0$. 
\end{proof}

Now we construct certain central extensions of $\pi_1 (X , x_0)$. With notation as in the preceeding section set $d_s = r$ and $d_{s_0} = 1$ and write $\Gamma_r$ for $\Gamma$. For $g = 0$ we have $\Gamma_r = \oN = 1$. For $g \ge 1$ there is the central extension \eqref{eq:8} where we have identified $\oN = \overline{I}_s$ with $\mu_r$
\begin{equation}
  \label{eq:10}
  1 \longrightarrow \mu_r \longrightarrow \Gamma_r \longrightarrow \pi_1 (X , x_0) \longrightarrow 1 \; .
\end{equation}
One can show that this extension corresponds to $1 \in \Z / r$ under the canonical isomorphisms
\[
H^2 (\pi_1 (X , x_0) , \mu_r) = H^2 (X , \mu_r) = \Z / r \; .
\]

We elaborate on this a bit:

{\bf Digression} (without proofs) More generally let $S \subset X$ be a finite set of closed points and set $U = X \ohne S$. Assuming $g \ge 1$ let $d_2$ be the transgressions in the Hochschild--Serre spectral sequence for the extension
\begin{equation}
  \label{eq:11}
  1 \longrightarrow N\longrightarrow \pi_1 (U ,x_0) \longrightarrow \pi_1 (X,x_0) \longrightarrow 1
\end{equation}
and in the Leray spectral sequence for the inclusion $\lambda : U \hookrightarrow X$. Then we have a canonical commutative diagram
\[
\xymatrix{
H^0 (\pi_1 (X , x_0) , H^1 (N , \mu_r)) \ar[r]^-{d_2} \ar@{=}[d] & H^2 (\pi_1 (X ,x_0) , \mu_r) \ar@{=}[d]\\
H^0 (X , R^1 \lambda_* \mu_r) \ar[r]^-{d_2} \ar[d]^{\wr} & H^2 (X , \mu_r) \ar@{=}[d] \\
H^2_S (X , \mu_r) \ar[r] & H^2 (X , \mu_r) \; .
} 
\]
Using the identifications
\[
H^0 (\pi_1 (X,x_0) , H^1 (N , \mu_r)) = \Hom_{\pi_1 (X, x_0)} (N^{\abb} , \mu_r) = \Hom (\tN , \mu_r)
\]
where $\tN$ was defined before theorem \ref{t9} we get a commutative diagram
\[
\xymatrix{
\Hom (\tN , \mu_r) \ar[r]^-{d_2} \ar[d]_{\omega} & H^2 (\pi_1 (X ,x_0) , \mu_r) \ar@{=}[d] \\
\bigoplus\limits_{s \in S} \Z / r \ar[r]^{\Sigma} & \Z / r \; .
}
\]
Here for a homomorphism $\varphi : \tN \to \mu_r$ the element $d_2 (\varphi)$ is the class of the push-out of \eqref{eq:11} via $\varphi$. The map $\omega$ is given as follows. Identifying the inertia groups at the points $s \in S$ with $\mu_r$ we have $\tN = \prod_{S \in S} \mu_r$. Then $\omega$ is the composition:
\[
\Hom (\tN , \mu_r) = \bigoplus_{s \in S} \Hom (\mu_r , \mu_r) = \bigoplus_{s \in S} \Z / r \; .
\]
The extension \eqref{eq:10} for example is the push-out of \eqref{eq:11} via the map $\varphi \in \Hom (\tN , \mu_r)$ corresponding to $\omega (\varphi) = (1,0)$ Hence it corresponds to $\Sigma (1,0) = 1 \in \Z / r$.

Let $E$ be a vector bundle in $\eB_X^\mu$, i.e. $E$ is a vector bundle on $X_{\C_p}$ of slope $\mu$ with potentially strongly semistable reduction in the sense of definition \ref{t2}. We denote by $r$ the rank of $E$ and by $d$ its degree, so that $\mu =d / r$. 
Consider a finite Galois covering $\alpha: Y \rightarrow X$ with group $G = \mu_r$ as in construction \ref{t11}.
We denote the unique points in $Y$ lying over ${s}$ and ${s_0}$ by $t$ and $t_0$, respectively. Let $\Oh(-d {t})$ be the line bundle on $Y$ associated to the divisor $-d {t}$. For all $\sigma \in G$ there is a canonical isomorphism of line bundles $\sigma^\ast \Oh(-d {t}) \rightarrow \Oh(-d {t})$. Hence the line bundle $\Oh(-d {t})$
 carries a natural $G$-action. The group $G_{{t_0}}$ acts trivially on the fibre of $\Oh(-d {t})$ over ${t_0}$ and $G_{{t}} = G$ acts on the fibre of $\Oh (-dt)$ over $t$ by a character $\chi : G \to \mu_r$. 

The line bundle $\Oh(-d {t})$ has slope $-d$. Since $\alpha$ has degree $r$, the pullback bundle $\alpha^\ast_{\C_p} E$ has slope $d$. By corollary \ref{t5}, the bundle $F = \alpha^\ast_{\C_p} E \otimes \Oh(-d {t})$ is contained in the category $\eB_Y^0$. It is naturally equipped with isomorphisms $a_{\sigma} : \sigma^* F \to F$ for $\sigma \in G$ satisfying the cocycle relations. The group $G_t$ acts by the character $\chi$ on $F_t$ and $G_{t_0}$ acts trivially on $F_{t_0}$. By theorem \ref{t9} and remark \ref{t10n} the representation $\rho_F$ of $\Pi_1 (V)$ induces a continuous representation
\[
\rho : \Gamma_r \to \GL (\rho (x_0))
\]
of the group $\Gamma_r$ from sequence \eqref{eq:10} above. 

We have by definition:
\[
\rho (x_0) = \{ (f_y) \in \prod_{\alpha (y) = x_0} (\alpha^* E)_y \otimes \Oh (-dt)_y \tei a_{\sigma,y} (f_{\sigma y}) = f_y \quad \mbox{for all} \; \sigma \in G \} \; .
\]
Since $\Oh (-dt)$ restricted to $V$ is the pullback of $\Oh_U$, we can identify $\rho (x_0)$ with $E_{x_0}$. Hence from $\rho$ we obtain the desired representation 
\[
\rho_{E, x_0} : \Gamma_r \to \GL (E_{x_0}) \; .
\]

\begin{prop}
  \label{t12n}
Assuming $g \ge 1$ there is some $a \in (\Z / r)^*$ such that we have
\[
\rho_{E,x_0} (\zeta) = \zeta^{ad} \id_{E_{x_0}} \quad \mbox{for all} \; \zeta \in \mu_r \subset \Gamma_r \; .
\]
\end{prop}

\begin{proof}
  By theorem \ref{t9} and remark \ref{t10n} there is a generator $\tau$ of $G_t = G$ such that we have $\rho_{E,  x_0} (\overline{\gamma}_s) = \chi (\tau)$. On the other hand $a_{\tau , t}$ acts on the fibre of $\Oh (-dt) = \Oh (-t)^{\otimes d}$ over $t$ by multiplication with $\eta^d$ where $\eta$ is a primitive $r$'th root of unity. Hence we have $\chi (\tau) = \eta^d$. This implies the assertion.
\end{proof}

\end{document}